\documentclass[10pt]{article}

\usepackage{amsthm}
\usepackage{amsfonts}
\usepackage{a4wide}
\newtheorem{theorem}{Theorem}
\newtheorem{lemma}[theorem]{Lemma}
\newtheorem{corollary}[theorem]{Corollary}
\newtheorem{claim}{Claim}

\begin{document}
\title{Algebraic proof of Brooks' theorem\thanks{This research was partially supported by the grant GA {\v C}R 201/09/1097.}}
\author{Jan Hladk{\'y}\thanks{Zentrum Mathematik, Technische Universit\"at M\"unchen, Boltzmannstra\ss e 3, D-85747 Garching bei M\"unchen, Germany, and Department of Applied Mathematics, Faculty of Mathematics and Physics, Charles University, Malostransk\'e n\'am\v{e}st\'{\i}~25, 118~00 Prague, Czech Republic. E-mail: {\tt hladk@seznam.cz}.} \and
        Daniel Kr{\'a}l'\thanks{Institute for Theoretical Computer Science, Faculty of Mathematics and Physics, Charles University, Malostransk\'e n\'am\v{e}st\'{\i}~25, 118~00 Prague, Czech Republic. E-mail: {\tt kral@kam.mff.cuni.cz}. The Institute for Theoretical Computer Science is supported by Ministry of Education of the Czech Republic as project 1M0545.}\and
	Uwe Schauz\thanks{Department of Mathematics and Statistics, King Fahd University of Petroleum and Minerals, Dhahran 31261, Saudi Arabia. E-mail: {\tt schauz@kfupm.edu.sa}}}
\date{}
\maketitle
\begin{abstract}
We give a proof of Brooks' theorem and its list coloring extension
using the algebraic method of Alon and Tarsi; this also shows that
the Brooks' theorem remains valid in a more general game coloring setting.
\end{abstract}

\section{Introduction}

One of the most famous theorems on graph colorings is Brooks' theorem~\cite{bib-brooks41}
which asserts that every connected graph $G$ with maximum degree $\Delta$ is $\Delta$-colorable
unless $G$ is an odd cycle or a complete graph. Brooks' theorem has been extended in various
directions, e.g., its list version can be found in~\cite{bib-vizing76}, also see~\cite{bib-kostochka96+}.
A short and elegant proof of Brooks' theorem was given in~\cite{bib-lovasz75} by Lov{\'a}sz and
an algebra-based proof in~\cite{bib-tverberg} by Tverberg.
In this note, we give yet another proof of Brooks' theorem.
Our proof is based on the algebraic method of Alon and Tarsi
from \cite{bib-alon92+} which also yields its list coloring version. In fact,
the argument shows that Brooks' theorem remains valid in a more general
setting of a coloring game of two players from~\cite{bib-uwe09,bib-uwe}, which we later describe.
Our proof can also be seen as the first step towards obtaining a variant of Brooks' theorem
for circular colorings.

The algebraic method of Alon and Tarsi is based on studying of properties of a certain graph polynomial.
This polynomial is closely related to the existence of special orientations
with bounded in-degrees. We summarize this relation in the next theorem.
\begin{theorem}[Alon and Tarsi~\cite{bib-alon92+}]
\label{thm-alon-tarsi}
Let $G$ be an oriented graph and $d^+_v$ the indegree of a vertex $v$.
If the numbers of even and odd Eulerian subgraphs of $G$ differ,
then $G$ can be colored from any lists $L(v)$ such that $|L(v)|\ge d^+_v+1$.
In particular, if the maximum in-degree of $G$ is $k$, $G$ is list $(k+1)$-colorable.
\end{theorem}
Let us remind that an {\em even Eulerian subgraph} is a spanning subgraph of $G$ with even number of edges
such that each vertex has the same in-degree and out-degree.
Similarly, an odd Eulerian subgraph is such a subgraph
with an odd number of edges. Even and odd Eulerian subgraphs do not need to be connected and can contain
isolated vertices. Theorem~\ref{thm-alon-tarsi} has been successfully applied to several coloring problems,
e.g., showing that planar bipartite graphs are list $3$-colorable~\cite{bib-alon92+}.
We refer the reader to~\cite{bib-alon99} and \cite{bib-uwe08} for further applications
of the algebraic tool behind this theorem.

In~\cite{bib-uwe}, the third author showed that the assumption of Theorem~\ref{thm-alon-tarsi}
implies a stronger graph property than list $(k+1)$-colorability, namely, $(k+1)$-paintability.
Consider the following game played by two players, Mr.~Paint and Mrs.~Correct. At the beginning,
each vertex of a graph $G$ is equipped with $k$ erasers. The players move in turns, Mr.~Paint starting.
In his turn, Mr.~Paint marks a non-empty subset of vertices; Mrs.~Correct then removes an independent
subset of marked vertices from the graph and clears marks of other vertices with their erasers.
Mr.~Paint wins if Mrs.~Correct is unable to clear the mark of an unremoved vertex because all
its erasers has been previously used. Mrs.~Correct wins if all the vertices are eventually removed
from the graph. If Mrs.~Correct always wins, the graph is said to be $(k+1)$-paintable.
The following theorem can be found in~\cite{bib-uwe}:
\begin{theorem}
\label{thm-paint}
Let $G$ be an oriented graph and $d^+_v$ the indegree of a vertex $v$.
If the numbers of even and odd Eulerian subgraphs of $G$ differ,
then Mrs.~Correct always wins the game if every vertex $v$
is initially equipped with $d^+_v$ erasers.
In particular, if the maximum in-degree of $G$ is $k$, $G$ is $(k+1)$-paintable.
\end{theorem}
It is not hard to observe that every $(k+1)$-paintable graph is also list $(k+1)$-colorable.
The proof of Theorem~\ref{thm-paint} given in~\cite{bib-uwe} is purely combinatorial and
thus it gives a combinatorial proof of Theorem~\ref{thm-alon-tarsi}.

Another motivation for an algebraic proof of Brooks' theorem
stems from the area of circular colorings introduced in~\cite{bib-vince88}.
Circular colorings attracted a considerable amount of interest of researchers (see
two recent surveys~\cite{bib-zhu01,bib-zhu06} on the topic by Zhu). Classical and circular colorings
are closely related, in particular, it holds that $\chi(G)=\lceil\chi_c(G)\rceil$ (we omit a formal
definition of the circular chromatic number as it is not needed for further exposition).
An analogous equality
is not true for their list counterparts. The circular list chromatic number is always at least the list
chromatic number decreased by one but it is not upper-bounded by any function of the list chromatic
number~\cite{bib-zhu05}.

Circular list colorings seem to be of surprising difficulty, e.g., Norine~\cite{bib-norine08} only recently
proved that the list chromatic number of even cycles is equal to two. In his proof, he has successfully
applied the algebraic method of Alon and Tarsi from~\cite{bib-alon92+}. For another application of this
method to circular list colorings, see~\cite{bib-norine08+,bib-norine08++}. It seems natural to ask whether this approach
can also be used to prove the variant of Brooks' theorem for circular list colorings, which is still
not known~\cite{bib-zhu05}, and the natural
first step towards this goal is finding a proof of the classical Brooks' theorem based on the Alon and Tarsi method
which we present in this short note.

\section{Structural lemma}

In this section, we give a structural lemma which allows us to apply the algebraic technique of Alon and Tarsi
in our proof of Brooks' theorem.
Recall that a Gallai tree is a graph whose every block is a complete graph or an odd cycle.

\begin{lemma}
\label{lm-struct}
Every connected graph $G$ that is not a Gallai tree contains an even cycle $C$ with at most one chord
as an induced subgraph.
\end{lemma}

\begin{proof}
We prove the lemma in a series of five claims.
\begin{claim}
\label{cl-1}
The graph $G$ contains a cycle $C$ that induces neither a complete graph nor an odd cycle.
\end{claim}
Let $H$ be a block of the block-decomposition of $G$ that
is neither a complete graph nor an odd cycle (if $G$ is $2$-connnected, set $H$ to be $G$).

Let $C_0$ be the longest cycle of $H$.
Suppose that $C_0$ induces a complete graph or an odd cycle.
Since $H$ is neither a complete graph nor an odd cycle,
$C_0$ does not contain all vertices of $H$.
Hence, there exists a vertex $v$ of $C_0$ with a neighbor $v'$ not on $C_0$.
Let $P$ be the shortest path from $v'$ to $C_0\setminus v$ in $H\setminus v$
which exists because $H$ is $2$-connected. Let $w$ be the end of $P$ on $C_0$.

If $C_0$ induced a complete graph, we could assume that
$v$ and $w$ are neighbors on $C_0$. However, if $v$ and $w$ are neighbors,
then the cycle $C_0$ could be prolonged by removing the edge $vw$ and
adding the edge $vv'$ and the path $P$
which would contradict the choice of $C_0$.
We conclude that $C_0$ induces an odd cycle and
the vertices $v$ and $w$ are not neighbors on $C_0$.

Let $C_1$ and $C_2$ be the two cycles formed by the edge $vv'$, the path $P$ and
one of the two parts of $C_0$ delimtied by $v$ and $w$. Since $C_0$ is odd,
one of the cycles $C_1$ and $C_2$ is even, say $C_1$. The cycle $C_1$ does not
induce a complete graph in $H$ because the vertices $v$ and $w$ are not adjacent and
thus $C_1$ is the sought cycle.

\begin{claim}
\label{cl-2}
Let $C=v_1\ldots v_{\ell}$ be the shortest cycle with the properties given in Claim~\ref{cl-1}.
If $\ell\ge 5$, then no vertex of $C$ is adjacent to all other vertices on $C$.
\end{claim}

Suppose that the vertex $v_1$ is adjacent to all the vertices $v_2,\ldots,v_{\ell}$.
By the choice of $C$, both the cycles $v_1\ldots v_{\ell-1}$ and $v_1v_3\ldots v_{\ell}$
induce complete graphs. By the choice of $C$, the vertices $v_2$ and $v_{\ell}$ are not adjacent.
However, the cycle $v_1v_2v_3v_{\ell}$ is then a shorter cycle satisfying the properties of Claim~\ref{cl-1}.

\begin{claim}
\label{cl-3}
Let $C=v_1\ldots v_{\ell}$ be the shortest cycle with the properties given in Claim~\ref{cl-1}.
Each vertex of $C$ is incident with at most one chord.
\end{claim}

If $\ell=4$, then $C$ is a cycle of length four with at most one chord.
Assume that $\ell\ge 5$ and that there exists a vertex, say $v_1$, adjacent
to vertices $v_a$ and $v_b$, $3\le a<b\le\ell-2$.
However,
the cycle $v_1\ldots v_b$ or the cycle $v_1v_a\ldots v_{\ell}$ does not
induce a complete graph (since $v_1$ is not adjacent to all vertices on $C$ by Claim~\ref{cl-2})
or a cycle (since it has a chord $v_1v_a$ or $v_1v_b$). Hence, it is a shorter cycle
with the properties given in Claim~\ref{cl-1} which contradicts the choice of $C$.

\begin{claim}
\label{cl-3+}
The shortest cycle $C=v_1\ldots v_{\ell}$ with the properties given in Claim~\ref{cl-1} is even.
\end{claim}

Suppose that $C$ is odd. Hence, $C$ must have a chord, say $v_1v_i$.
One of the cycles $v_1\ldots v_i$ and $v_1v_i\ldots v_{\ell}$ is even;
by symmetry we can assume that $v_1\ldots v_i$ is even. By the choice
of $C$, $v_1,\ldots,v_i$ must induce a complete graph and thus $v_1$
is adjacent to $v_3$ and $v_i$ which is impossible by Claim~\ref{cl-3}.

\begin{claim}
\label{cl-4}
The shortest cycle $C=v_1\ldots v_{\ell}$ with the properties given in Claim~\ref{cl-1}
contains at most one chord.
\end{claim}
Suppose that $C$ has at least two chords, say $v_av_b$ and $v_cv_d$.
If the two chords do not cross, we can assume without loss of
generality that $1\le a<b<c<d\le\ell$. But then the cycle $v_1\ldots v_av_b\ldots v_{\ell}$
would be a shorter cycle with the properties given in Claim~\ref{cl-1}. 
Hence, the chords $v_av_b$ and $v_cv_d$ cross.
By symmetry, we can assume that $1\le a<c<b<d\le\ell$.

The vertices $v_a$, $v_b$, $v_c$ and $v_d$ split the cycle $C$ into four parts;
let $n_{xy}$, $x\in\{a,b\}$ and $y\in\{c,d\}$, be the number of vertices of $C$
between $v_x$ and $v_y$. If not all $n_{xy}$ have the same parity, then there are two
consecutive parts (viewed in the order they correspond to the parts of $C$) with different
parities, say the parities of $n_{ac}$ and $n_{bc}$ are different. Then,
the cycle $v_a\ldots v_c\ldots v_b$ is an even cycle satisfying the properties of Claim~\ref{cl-1}
which is shorter than $C$. Since this is impossible, we can assume that the parities of
all the four numbers $n_{xy}$ are the same and one, say $n_{ad}$, is non-zero.
Now, the shorter even cycle $v_a\ldots v_cv_d\ldots v_b$ must induce a complete
graph violating Claim~\ref{cl-3} for $C$.

\smallskip

Claims~\ref{cl-3+} and~\ref{cl-4} now imply the lemma.
\end{proof}

\section{Main result}

Before presenting our proof of Brooks' theorem, we need to recall a simple
folklore structural result on ordering vertices of a connected graph. We include its
proof for completeness.

\begin{lemma}
\label{lm-dfs}
Let $G$ be a connected graph and $v$ an arbitrary vertex of $G$. The vertices of $G$
can be ordered in such a way that every vertex except for $v$ is preceded by at least
one of its neighbors.
\end{lemma}

\begin{proof}
Consider an arbitrary spanning tree $T$ of $G$ and root it at $v$. The vertices of $G$
are ordered in the following way: the first vertex is the root $v$ followed by all
its children (vertices of the second level of $T$). Then, all vertices of the third level
are listed, then all vertices of the forth level, etc. Since each vertex except for $v$
is preceded by its parent, the obtained ordering has the desired property.
\end{proof}

We prove our main result in a more general setting of degree list colorability
considered in~\cite{bib-borodin77,bib-erdos79+} and degree paintability.
The list coloring version of Brooks' theorem asserts that every connected graph except
for Gallai trees can be colored from lists such that each vertex has a list of size equal
to its degree. If a graph can be colored from all such lists, we say that it is {\em list
degree colorable} (in other words, a connected graph is list degree colorable if and only
if it is not a Gallai tree). Similarly, we speak of a graph being degree paintable:
a graph is {\em degree paintable} if Mrs.~Correct always wins if the number of erasers
at each vertex is initially the same as its degree decreased by one.

We are now ready to give the algebraic proof of Brooks' theorem.

\begin{theorem}
\label{thm-main}
Let $G$ be a connected graph. If $G$ is not a Gallai tree,
then $G$ is list degree colorable and degree paintable.
\end{theorem}

\begin{proof}
It is enough to prove that $G$ has an orientation such that
each vertex has at least one out-going edge and
the numbers of even and odd Eulerian subgraphs differ.
Theorems~\ref{thm-alon-tarsi} and \ref{thm-paint}
would then imply the statement.

By Lemma~\ref{lm-struct},
$G$ contains an even cycle $C$ with at most one chord.
Contract $C$ to a single vertex $w$ in $G$ and apply Lemma~\ref{lm-dfs} to the resulting graph.
We obtain an ordering $v_1,\ldots,v_n$ of the vertices of $G$ by replacing $w$ in the order
given by Lemma~\ref{lm-dfs} with the vertices of $C$, inserted in an arbitrary order.
The edges of $C$ are oriented in a cyclic way.
Every edge $v_iv_j$, $i<j$,
that is not contained in $C$, is now oriented from $v_j$ to $v_i$ (this rule also applies
to the chord of $C$ if it exists).

There are always two even Eulerian subgraphs of the constructed orientation of $G$,
in particular, the empty one and the one formed by the cycle $C$. If $C$ contains a chord $e$,
there is also an Eulerian subgraph formed by the chord and one of the two parts delimited by $e$,
which is either even or odd. Since there are no other Eulerian subgraphs,
the numbers of odd and even Eulerian subgraphs must differ. Since every vertex has at least
one out-going edge, the statement of the theorem now follows.
\end{proof}

Theorem~\ref{thm-main} now yields the following:

\begin{corollary}
\label{cor-main}
Let $G$ be a connected graph with maximum degree $\Delta$. If $G$ is neither a complete graph nor
an odd cycle, then $G$ is list $\Delta$-colorable and $\Delta$-paintable.
\end{corollary}

\section*{Acknowledgement}

The third author gratefully acknowledges the
support provided by the King Fahd University of Petroleum and
Minerals during this research.

\end{document}